\documentclass[aps,pre,twocolumn,groupedaddress,amsmath,amssymb,floatfix]{revtex4-1}

\bibstyle{apsrev4-1}
\usepackage{graphicx}

\begin{document}

\title{Correlation transfer in stochastically driven oscillators over long and short time scales}

\author{Aushra Abouzeid and Bard Ermentrout}
\affiliation{University of Pittsburgh}

\date{\today}

\begin{abstract}
In the absence of synaptic coupling, two or more neural oscillators may become synchronized by virtue of the statistical correlations in their noisy input streams. Recent work has shown that the degree of correlation transfer from input currents to output spikes depends not only on intrinsic oscillator dynamics, but also depends on the length of the observation window over which the correlation is calculated. In this paper we use stochastic phase reduction and regular perturbations to derive the correlation of the total phase elapsed over long time scales, a quantity which provides a convenient proxy for the spike count correlation. Over short time scales, we derive the spike count correlation directly using straightforward probabilistic reasoning applied to the density of the phase difference. Our approximations show that output correlation scales with the autocorrelation of the phase resetting curve over long time scales. We also find a concise expression for the influence of the shape of the phase resetting curve on the initial slope of the output correlation over short time scales. These analytic results together with numerical simulations provide new intuitions for the recent counterintuitive finding that type I oscillators transfer correlations more faithfully than do type II over long time scales, while the reverse holds true for the better understood case of short time scales.
\end{abstract}

\maketitle

While the jury is still out on the functional role of synchrony and correlations in neural firing, the ubiquity of these phenomena in the nervous system is suggestive. One long-standing hypothesis holds that correlated activity in the visual system underlies feature binding. Synchronous oscillations may also play a role in amplifying signals \cite{tiesinga:2004}, transmitting information from one layer to another \cite{salinas_impact_2000, kuhn_higher-order_2003, tetzlaff_dependence_2008}, or such oscillations may encode information directly \cite{decharms_primary_1996, samonds_cooperation_2003, kohn_stimulus_2005, brent:nature, gray_oscillatory_1989, biederlack_brightness_2006, chacron_population_2008, josi_stimulus-dependent_2009}. On the other hand, correlations may negatively impact the signal-to-noise ratio \cite{zohary_correlated_1994, johnson_sensory_1980, britten_analysis_1992, bair_correlated_2001}, and excessive synchrony is a hallmark of neurological disorders such as epilepsy and Parkinson's disease. 

To understand the function of oscillatory correlations, or one day achieve clinically relevant control over them, we must first understand the underlying biophysical mechanisms. While synchrony can arise as the result of anatomical connectivity between neurons, much recent work \cite{galan:2006, goldobin:2005, teramae:2004, nakao:2005, ermentrout, galanfp} has brought to light ways in which correlated activity develops from the inherent stochastisicity of neural systems. Thus, in the absence of direct coupling, two or more neural oscillators may become synchronized by virtue of the statistical correlations in their noisy input streams -- a phenomenon we will refer to as stochastic synchrony.

For our analysis of stochastic synchrony, we appeal to the theory of weak coupling, which holds in the stochastic context provided the amplitude of the noise is sufficiently small. In particular, a number of groups \cite{teramae:2004, goldobin:2005, nakao:2005, yoshi} have proved that the phase reduction technique \cite{kuramoto:1984} can be applied to oscillators receiving additive noise. Thus, we reduce a noisily driven oscillator to a scalar differential equation describing the evolution of the phase. This so-called phase equation depends only on the properties of the noise and the oscillator's phase resetting curve (PRC) which characterizes how small perturbations influence the oscillator's subsequent timing or phase.

\begin{figure}
\begin{center}
\includegraphics[width=3in]{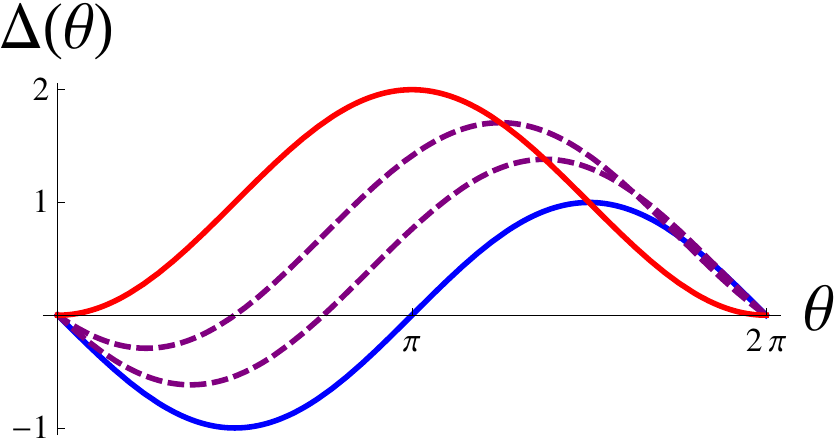}
\caption{We use the parametrization $\Delta(\theta)=-\sin(\theta+\alpha)+\sin(\alpha)$ to vary the PRC smoothly from type I (red), where $\alpha=\frac{\pi}{2}$ and $\Delta(\theta)=1-\cos(\theta)$, to type II (blue), where $\alpha=0$ and $\Delta(\theta)=-\sin(\theta)$. Note that intermediate values of $\alpha$ produce PRC shapes (dashed purple) that more closely resemble those found empirically \textit{in vivo}.}
\label{fig:PRCs}
\end{center}
\end{figure}

Neural oscillators can be classified into two types according to the bifurcations that occur as the dynamical system goes from a stable rest state to a stable limit cycle. Furthermore, the oscillator's bifurcation class has been shown to determine the shape of it's PRC and therefore it's ability to synchronize. Type I oscillators undergo the saddle-node-on-an-invariant-circle, or SNIC, bifurcation and the resulting PRC is strictly positive, indicating that perturbations can only advance the oscillator's phase. Type II cells undergo the Andronov-Hopf bifurcation, which produces a PRC with both negative and positive regions; typically, inputs occurring early in the cycle can delay the phase while later inputs advance it. See Fig.(\ref{fig:PRCs}).

An expanding body of work has demonstrated that over short time scales of less than one period, type II oscillators are more susceptible to stochastic synchrony than type I. This has been shown via simulations and \textit{in vivo} \cite{galan:2006, galan-optimal}, by deriving the probability distribution of the phase difference \cite{sashi}, by minimizing the Lyapunov exponent of the phase difference \cite{abouzeid:2009}, and most recently by calculating the spike count correlation over a range of time windows \cite{barreiro:2010}. The latter study further reports that this finding reverses over long timescales, namely that type I oscillators transmit correlations more faithfully than type II when observed over lengths of time much greater than one period.

In Section \ref{sec:1} we provide a brief introduction to the phase reduction technique in a stochastic setting. Next in Section \ref{sec:2} we use regular perturbations to give a novel and straightforward analysis of correlation transfer over long time scales. To facilitate our derivation, we use the total elapsed phase as a proxy for the spike count. Note that the total phase (modulo the period) and the spike count differ by at most one, which is a negligible quantity when many spikes have been observed over a long time window. The expression we derive for the correlation coefficient of the total phase agrees both qualitatively and quantitatively with the results found in \cite{barreiro:2010}.

In Section \ref{sec:3} we consider short time scales less than or equal to the period of the oscillation. In this case, the total phase cannot be used to approximate the spike count. We therefore derive the spike count correlation directly, using simple probabilistic reasoning applied to the density of the phase difference. Our analytic results together with Monte Carlo simulations corroborate earlier work showing type II oscillators transfer correlations more readily than type I over short time windows. 

%------------------------------------
%------------------------------------

%%%%%%%%%%%%%%%%%%%%%%%%%%%%%%%%%%%%%%%%%%%%%%%%%%%%
\section{Noisy oscillators}
\label{sec:1}
%%%%%%%%%%%%%%%%%%%%%%%%%%%%%%%%%%%%%%%%%%%%%%%%%%%%

Let us begin with a neural oscillator receiving additive noise with equations of motion given by
\begin{equation*}
dX = F(X) dt + \sigma \xi,
\end{equation*}
where $X\in \mathbb{R}^n$ and $\xi$ is a white noise process. When $\sigma=0$, we assume the noiseless system has an asymptotically stable periodic solution $X_0(t) = X_0(t+\tau)$ with period $\tau$.

As in the deterministic case, we can reduce this high-dimensional system to a scalar equation for the evolution of the phase $\theta$ around the limit cycle. Let $\phi :\mathbb{R}^n\rightarrow \mathbb{S}^1$ map a neighborhood of the limit cycle to the phase on a circle. That is, $\theta = \phi(X)$, with $\theta\in [0,1)$. Then $\theta$ satisfies
\begin{equation*}
\frac{d\theta}{dt} = 1 + \sigma \nabla_X\phi(X) \cdot \xi,
\end{equation*}
where we have normalized the unperturbed period to be one. Next we can close the equation by assuming the noise amplitude $\sigma$ is sufficiently small, so that the system trajectory can be approximated by the noiseless limit cycle $X_0$:
\begin{equation}
\dot{\theta} \approx 1 + \sigma Z(\theta) \cdot \xi,
\label{eq:reduction}
\end{equation}
where $Z(\theta) = \nabla_X\phi(X_0(\theta))$ is the adjoint, or phase-dependent sensitivity of the trajectory to perturbation along the limit cycle. In the case of a neural oscillator, we assume the noisy perturbations arise as the result of stochastic synaptic input, which influences only the voltage variable. Hence $Z(\theta)$ has only one nonzero component, which is proportional to the phase resetting curve $\Delta(\theta)$.

Thus far, we have used the conventional change of variables to obtain Eq.(\ref{eq:reduction}), which therefore must be understood as a stochastic differential equation (SDE) in the Stratonovich sense. In order to eliminate the correlation between $\theta$ and $\xi$ we must use the It\^{o} change of variables, which will introduce an additional drift term:
\begin{equation*}
\dot{\theta} = 1 + \sigma \Delta(\theta) \xi + \frac{\sigma^2}{2}\Delta'(\theta) \Delta(\theta).
\end{equation*}
Here $'$ denotes differentiation with respect to $\theta$. For a detailed discussion of phase reduction in noisy oscillators see \cite{teramae:2009}.

%%%%%%%%%%%%%%%%%%%%%%%%%%%%%%%%%%%%%%%%%%%%%%%%%%%%
\section{Correlation transfer over long time scales}
\label{sec:2}
%%%%%%%%%%%%%%%%%%%%%%%%%%%%%%%%%%%%%%%%%%%%%%%%%%%%

We now consider the transfer of correlations over time scales much larger than the natural period of the oscillators. Given the level of correlation between the noisy inputs, we wish to know what level of correlation remains between the spike count of two oscillators after some time. For analytic convenience, however, we will use the total phase that has elapsed as a proxy for the spike count. Since these quantities differ by at most one, the discrepancy will be negligible for the large spike counts that accrue over long time scales. 

Our system will consist of two identical phase oscillators receiving weak, correlated, but not identical, additive white noise. Keeping only terms up to order $\sigma$, we have
\begin{align}
\dot{\theta_1} &= 1+\sigma \Delta(\theta_1)\xi_1(t) \nonumber \\ 
\dot{\theta_2} &= 1+\sigma \Delta(\theta_2)\xi_2(t).
\label{eq:twoOscillators}
\end{align} 
The noise takes the form
\begin{align}
\xi_1 &= \sqrt{c} \text{ } \xi_{C} + \sqrt{1-c}  \text{ } \xi_A \nonumber \\
\xi_2 &= \sqrt{c} \text{ } \xi_{C} + \sqrt{1-c}  \text{ } \xi_B,
\label{eq:twoNoises}
\end{align}
where $\xi_A$, $\xi_B$ and $\xi_C$ are mutually independent, zero mean white noise processes, and $c\in [0,1]$ is the correlation between $\xi_1$ and $\xi_2$, which we will refer to as the input correlation.

%------------------------------------
\begin{figure*}
%\begin{center}
\includegraphics[width=6.5in]{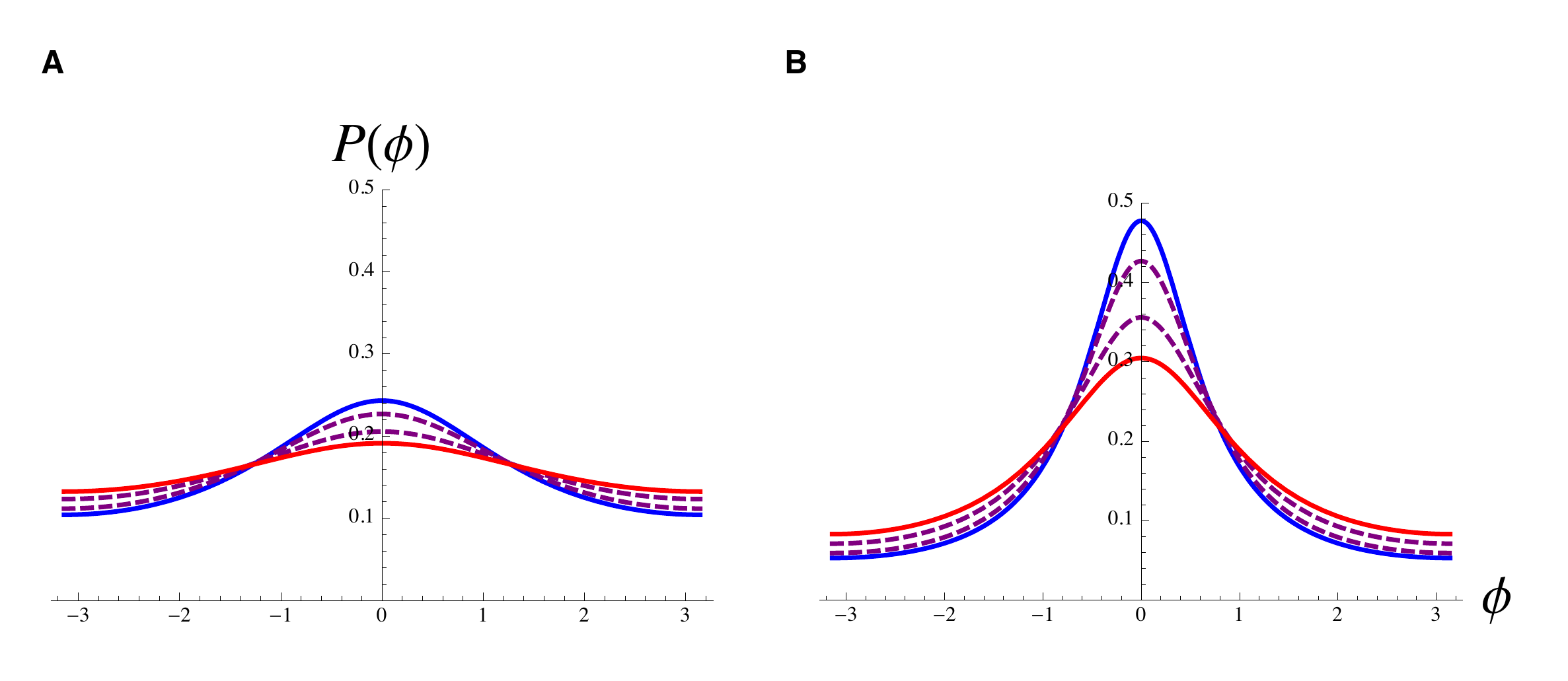}
\caption{The steady state distribution $P(\phi)$ of phase differences $\phi$ is shown for type I (red) and type II (blue) as well as for intermediate PRCs (dashed purple). Note that the unperturbed period of the oscillators is $2\pi$. (A) Input correlation $c=0.4$. (B) Input correlation $c=0.8$.}
\label{fig:probability}
%\end{center}
\end{figure*}
%------------------------------------

Next let us rewrite Eq.(\ref{eq:twoOscillators}) in the form of integral equations: 
\begin{align*}
\theta_1(t) &= t+ \theta_1(0)+\sigma \int_0^t \Delta(\theta_1(s))\xi_1(s) ds \\
\theta_2(t) &= t+ \theta_2(0)+\sigma \int_0^t \Delta(\theta_2(s))\xi_2(s) ds.
\end{align*}

Let $T$ be length of the window of time over which we will observe the system. Throughout this discussion we will assume  that our system has reached equilibrium, and that time has been reparametrized so that our observation takes place on the interval $t\in [0,T]$. In order to quantify the total phase traversed during this time, we subtract the initial phases by defining $q_i(T) = \theta_i(T) - \theta_i(0)$ for $i=1,2$. Thus the total phase traversed over a time window of length $T$ is given by:
\begin{equation*}
q_i(T) = T+\sigma \int_0^T \Delta(\theta_i(s))\xi_i(s) ds.
\end{equation*}
with $q_i(0)=0$  for $i=1,2$. Finally, since we assume $\sigma$ is small, let us simplify the integrands by expanding the phase to lowest order: 

\begin{equation}
\theta_i(t) = t+ \theta_i(0) + \mathcal{O}(\sigma).
\label{eq:lowestOrder}
\end{equation}
Then we have $\Delta(\theta_i(s)) = \Delta(t+ \theta_i(0))$, and thus 
\begin{equation}
q_i(T) = T+\sigma \int_0^T \Delta(s+\theta_i(0))\xi_i(s) ds
\label{eq:totalPhase}
\end{equation}

When taking expectations of the quantities in Eq.(\ref{eq:totalPhase}), we must keep in mind that there are four stochastic variables over which averaging must take place. In particular, we must average over the white noise signals $\xi_1(t)$ and $\xi_2(t)$ and the initial conditions $\theta_1(0)$ and $\theta_2(0)$.

Assuming we begin observation after the system has reached equilibrium, we can take one of the initial conditions, say $\theta_1(0)$, to be distributed uniformly on the interval $[0,2\pi ]$. However, at equilibrium the phases obey the steady state probability distribution $P(\phi)$ derived in \cite{sashi} and \cite{nakao:2007}, which depends only on the phase difference $\phi(t) = \theta_2(t) - \theta_1(t)$. Therefore, the average of Eq.(\ref{eq:totalPhase}) is computed as

\begin{align}
&\text{E}[q_i(T)] = \text{E}\left[T+ \sigma \int_0^T \Delta(s+x)\xi_i(s) ds \right] \nonumber \nonumber\\
  &=  
\frac{1}{2 \pi} \int_0^{2\pi}\int_0^{2\pi} P(y-x) \times \nonumber\\
& \indent \indent \left[T+\sigma \int_0^T \Delta(s+x)\left<\xi_i(s)\right> ds \right] dx dy \nonumber \\
  &=  T + \frac{\sigma}{2 \pi} \int_0^{2\pi}\int_0^{2\pi} P(y-x) \times \nonumber\\
  & \indent \indent \int_0^T \Delta(\theta_i(s))\left<\xi_i(s)\right> ds dx dy \nonumber \\
  &=  T,
  \label{eq:meanPhase}
\end{align}
where $2\pi$ is the unperturbed period of the oscillators, $P(\phi)$ is the steady state probability distribution of the phase difference, and $x$ and $y$ represent  the initial conditions $\theta_1(0)$ and $\theta_2(0)$, respectively. The last line follows because the white noises have zero mean.
 
Our goal is to compute the correlation of the total phase traversed by the two oscillators:
\begin{equation}
\text{Cor}[q_1,q_2] = \frac{\text{Cov}[q_1,q_2]}{\sqrt{\text{Var}[q_1]\text{Var}[q_2]}}.
\label{eq:corrCoef}
 \end{equation}
 
\begin{widetext}
First, we derive the covariance as follows
 \begin{align*}
\text{Cov}[q_1,q_2](T) 
&= \text{E}[(q_1(T)-\text{E}[q_1(T)])(q_2(T)-\text{E}[q_2(T))]] \\
&= \text{E}[(q_1(T)-T)(q_2(T)-T)] \\
&= \text{E}\left[ \sigma^2 \int_0^T \Delta(s+\theta_1(0))\xi_1(s) ds  \int_0^T \Delta(s'+\theta_2(0))\xi_2(s') ds' \right] \\
&= \sigma^2\frac{1}{2\pi} \int_0^{2\pi} \int_0^{2\pi} P(y-x)\int_0^T \int_0^T \Delta(s+x)\Delta(s'+y)\left<\xi_1(s) \xi_2(s')\right> ds ds' dx dy \\
&= \sigma^2\frac{c_{in}}{2\pi} \int_0^{2\pi} \int_0^{2\pi} P(y-x)\int_0^T \int_0^T \Delta(s+x)\Delta(s'+y) \delta(s-s') ds ds' dx dy \\
&= \sigma^2\frac{c_{in}}{2 \pi} \int_0^{2\pi}\int_0^{2\pi} P(y-x)\int_0^T \Delta(s + x) \Delta(s + y) ds dx dy.
\end{align*}
\end{widetext}
Similarly, we find the variance to be
\begin{align*}
& \text{Var}[q_1](T) = \text{E}[(q_1(T)-\text{E}[q_1(T)])^2] \\
  & \indent = \sigma^2\frac{1}{2 \pi} \int_0^{2\pi}\int_0^{2\pi} P(y-x)\int_0^T \Delta(s + x)^2 ds dx dy.
\end{align*}
Note that we therefore have $\text{Var}[q_1] = \text{Var}[q_2]$, and hence the denominator of Eq.(\ref{eq:corrCoef}) can be simplified: $\sqrt{\text{Var}[q_1]\text{Var}[q_2]} = \text{Var}[q_1]$. This gives the total phase correlation as
\begin{align}
&\text{Cor}[q_1,q_2](T) \nonumber \\
& \indent= c  \text{ } \frac{\int_0^{2\pi}\int_0^{2\pi} P(y-x)\int_0^T \Delta(s + x) \Delta(s + y) ds dx dy}{\int_0^{2\pi}\int_0^{2\pi} P(y-x)\int_0^T \Delta(s + x)^2 ds dx dy}.
\label{eq:totalPhaseCorrelation}
\end{align}

%------------------------------------
\begin{figure}
%\begin{center}
\includegraphics[width=3.25in]{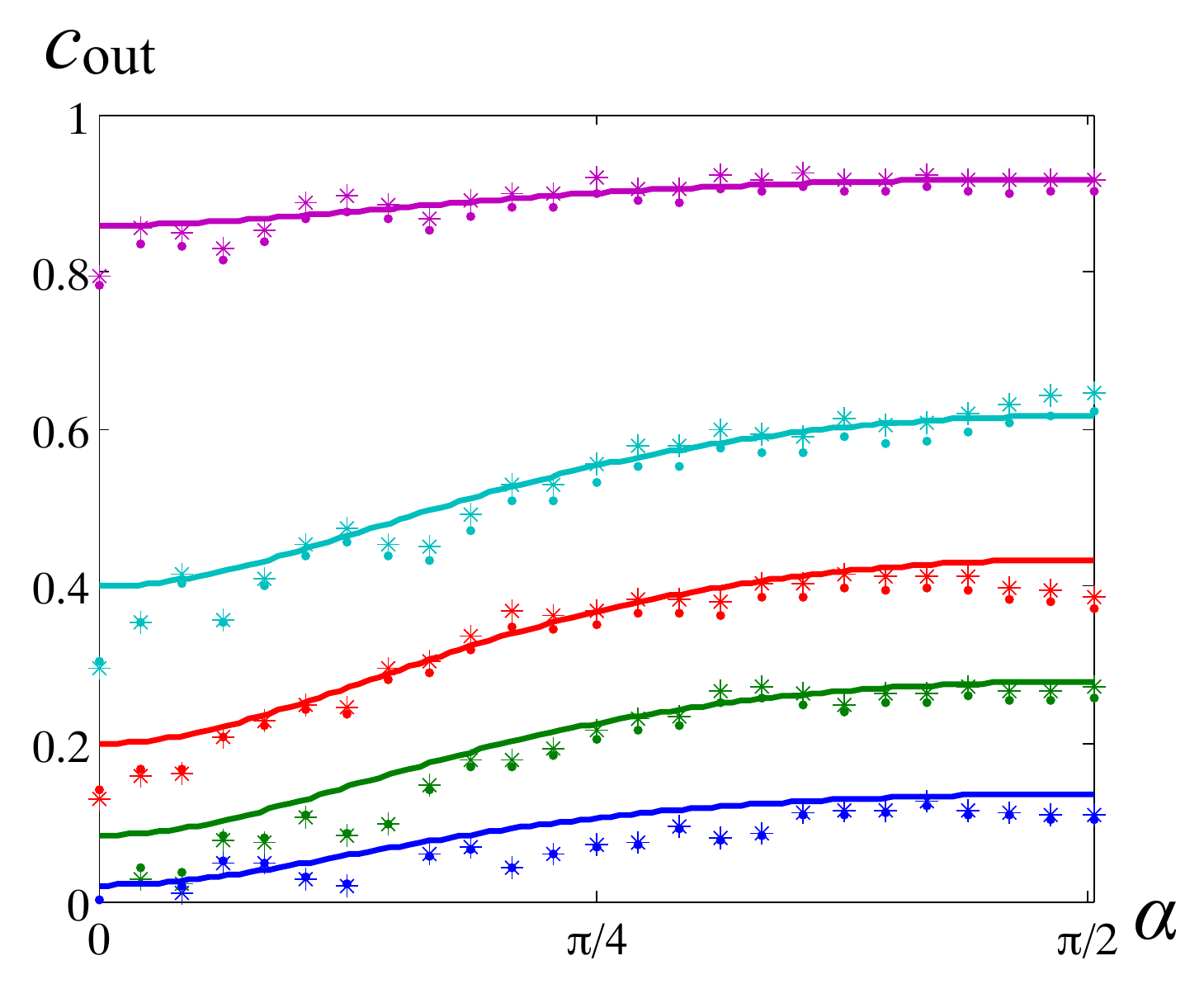}
\caption{Output correlation for large time windows is shown as a function of the PRC shape parameter $\alpha$. Note that when $\alpha=0$ the PRC is a pure sinusoid and therefore the oscillator is type II; when $\alpha=\pi/2$, the oscillator is type I (see Eq.(\ref{eq:paramPRC})). Theoretical curves (solid) are a good match for both the simulated total phase correlation (dotted) and the simulated spike count correlation (starred). Colors indicate the level of input correlation: 0.2 (blue), 0.4 (green), 0.6 (red), 0.8 (cyan), 0.99 (purple). In all cases, noise amplitude $\sigma = 0.05$.}
\label{fig:largeT}
%\end{center}
\end{figure}
%------------------------------------

Now let $h(x) = \int_0^{2\pi} \Delta (y) \Delta(y+x) dy$ be the autocorrelation of the PRC, and let $\phi(t) = \theta_2(t)-\theta_1(t)$ represent the phase difference as before. Then we can rewrite Eq.(\ref{eq:totalPhaseCorrelation}) as 

\begin{equation*}
c_{out} := \text{Cor}[q_1,q_2](T) = c  \text{ } \frac{\int_0^{2\pi} P(\phi) h(\phi) d\phi}{\int_0^{2\pi} P(\phi) h(0) d\phi}.
\end{equation*}
Note that the right hand side no longer depends on $T$ after we switched the order of integration and canceled the resulting factors of $T$ in both numerator and denominator. Next we can do away with the denominator entirely, since $h(0)$ does not depend on $\phi$, which leaves simply
\begin{equation}
c_{out} =\int_0^{2\pi} P(\phi) c \frac{h(\phi)}{h(0)} d\phi.
\label{eq:totalPhaseCorrelationSimplified}
\end{equation}

An expression for the steady-state probability density of the phase difference $P(x)$ was derived by Marella and Ermentrout in \cite{sashi}. Specifically, we have

\begin{equation*}
P(\phi) = \frac{N}{G(\phi)},
\end{equation*}
where $G(x) = 1-c \left( h(x)/h(0) \right)$, and $N$ is a normalizing constant, $N = 1/ \int_0^{2\pi} 1/G(x) dx$.
Let us further define the PRC to be
%------------------------------------
\begin{figure}
\begin{center}
\includegraphics[width=3.25in]{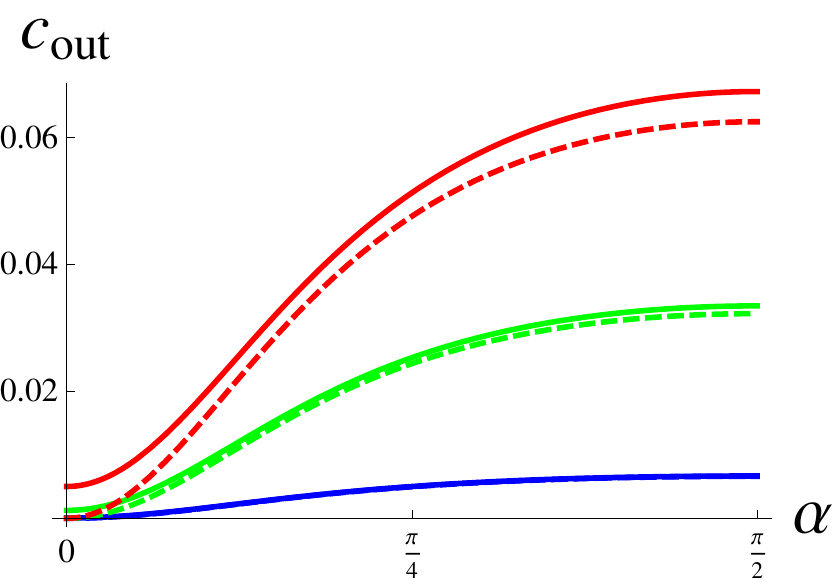}
\caption{The perturbation expansion of $c_{out}$ for small input correlation (dashed) agrees well with the full output correlation (solid). Note that, to lowest order in $c_{in}$, the output correlation goes to zero as the PRC shape parameter $\alpha$ goes to zero, that is, as the PRC shape approaches the pure type II. Colors indicate the level of input correlation: 0.01 (blue), 0.05 (green), 0.1 (red).}
\label{fig:approx}
\end{center}
\end{figure}
%------------------------------------
\begin{equation}
\Delta(\theta;\alpha) = -\sin(\theta + \alpha) - \sin(\alpha),
\label{eq:paramPRC}
\end{equation}
where $\alpha$ is a parameter that allows us to vary the PRC shape smoothly between type I ($\alpha=\pi/2$) and type II ($\alpha=0$). See Fig.(\ref{fig:PRCs}). Using this, the phase distribution over long time scales becomes a function of input correlation and the PRC shape parameter:

\begin{equation}
P(\phi;c,\alpha) = \frac{\sqrt{(c-1)(\cos(2\alpha) - 2)(2+ (c-1)\cos(2\alpha))}}{2\pi (2-c+(c-1)\cos(2\alpha) - c\cos(\phi))}.
\label{eq:ssPhaseDistribution}
\end{equation}

In the special cases where $\alpha=\pi/2$ and $\alpha=0$, Eq.(\ref{eq:paramPRC}) and Eq.(\ref{eq:ssPhaseDistribution}), together with Eq.(\ref{eq:totalPhaseCorrelation}), yield

%\begin{minipage}[b]{.9\linewidth}
%\setlength{\columnsep}{1pt}
%\begin{multicols}{2}
\begin{equation*}
\text{\underline{Type I}}
\end{equation*}
\begin{align}
\Delta_I(x) &= 1 - \cos(x) \nonumber\\
P_I(\phi;c) &= \frac{\sqrt{3}}{2 \pi} \frac{\sqrt{c^2 - 4 c + 3}}{(3 - 2 c - c \cos(\phi))}  \label{eq:PI} \\
c_{out,I} &= 1 - \frac{1}{3} \sqrt{3 (c-3) (c-1)} \nonumber
\end{align}
\begin{equation*}
\text{\underline{Type II}}
\end{equation*}
\begin{align}
\Delta_{II}(x) &= -\sin(x) \nonumber\\
P_{II}(\phi;c) &= \frac{1}{2 \pi} \frac{\sqrt{1 - c^2}}{(1 - c \cos(\phi))} \label{eq:PII}\\
c_{out,II} &= 1 - \sqrt{1 - c^2} \nonumber
 \end{align}
%\end{multicols}
%\end{minipage}
\vspace{0.5cm}

As in \cite{barreiro:2010}, we see in Fig.(\ref{fig:largeT}) that type I oscillators display greater output correlation than type II oscillators for any fixed value of the input correlation $c$, a surprising finding in light of earlier results that demonstrated the opposite relationship over short windows of observation.

Our intuition for this finding can be honed by performing a further perturbation expansion, now assuming small input correlation. For sufficiently small $c$, we can make the approximation

\begin{equation*}
\frac{1}{G(x)} = \frac{1}{1-c\frac{h(x)}{h(0)}} \approx 1+c\frac{h(x)}{h(0)}.
\end{equation*}
When we substitute this into Eq.(\ref{eq:totalPhaseCorrelationSimplified}) we find

\begin{equation}
c_{out} =  c \frac{\tilde{N}}{h(0)} \int_0^{2\pi} h(\phi) d\phi + \mathcal{O} (c ^2),
\label{eq:totalPhaseCorrelationLowest}
\end{equation}
where $\tilde{N} = 1/\int_0^{2\pi} \left(1+c_{in} h(x)/h(0) \right) dx$ is likewise approximated to lowest order in $c$. 

The form of Eq.(\ref{eq:totalPhaseCorrelationLowest}) demonstrates that output correlation scales with the integral of the PRC autocorrelation, and for the parametrized PRC in Eq.(\ref{eq:paramPRC}) we have

\begin{equation*}
\int_0^{2\pi} h(\phi) d\phi = 4 \pi^2 \sin(\alpha)^2.
\end{equation*}
In particular, $\alpha=0$  for the type II PRC, and hence $c_{out}=0$ to lowest order. Clearly, we have nonzero autocorrelation for nonzero $\alpha \leq \frac{\pi}{2}$, and hence PRCs that deviate from pure type II will produce higher output correlation over the long timescales considered here.

Expanding the remaining terms in Eq.(\ref{eq:totalPhaseCorrelationLowest}), we find the approximated output correlation takes the form
\begin{equation}
c_{out} =  \frac{2 c \sin(\alpha)^2}{2+c-(1+c) \cos(2 \alpha)}.
\end{equation}
In Fig.(\ref{fig:approx}) we see that this approximation agrees with Eq.(\ref{eq:totalPhaseCorrelation}) for $c=0.01$ and  $0.05$ but diverges for  $c=0.1$. Note that these curves would all lie below the lowest curve plotted in Fig.(\ref{fig:largeT}) if shown on the same scale.

We verify the preceding analysis by simulating two phase oscillators perturbed by additive white noise as described in Eq.(\ref{eq:twoOscillators}) and Eq.(\ref{eq:twoNoises}). The simulations used noise amplitude $\sigma = 0.05$, and the input correlation took the values $c\in \{0.2, 0.4, 0.6, 0.8, 0.99\}$. 

We computed the correlation coefficient of both the total phase and the spike count, using a range of observation windows $T$. As shown in Fig.(\ref{fig:largeT}), the total phase correlation and the spike count correlation agree closely both with each other and with the theoretical curves as a function of the PRC shape parameter $\alpha$.  

 %%%%%%%%%%%%%%%%%%%%%%%%%%%%%%%%%%%%%%%%%%%%%%%%%%%%
\section{Short time scales}
\label{sec:3}
 %%%%%%%%%%%%%%%%%%%%%%%%%%%%%%%%%%%%%%%%%%%%%%%%%%%%

%------------------------------------
\begin{figure}
\begin{center}
\includegraphics[width=3.5in]{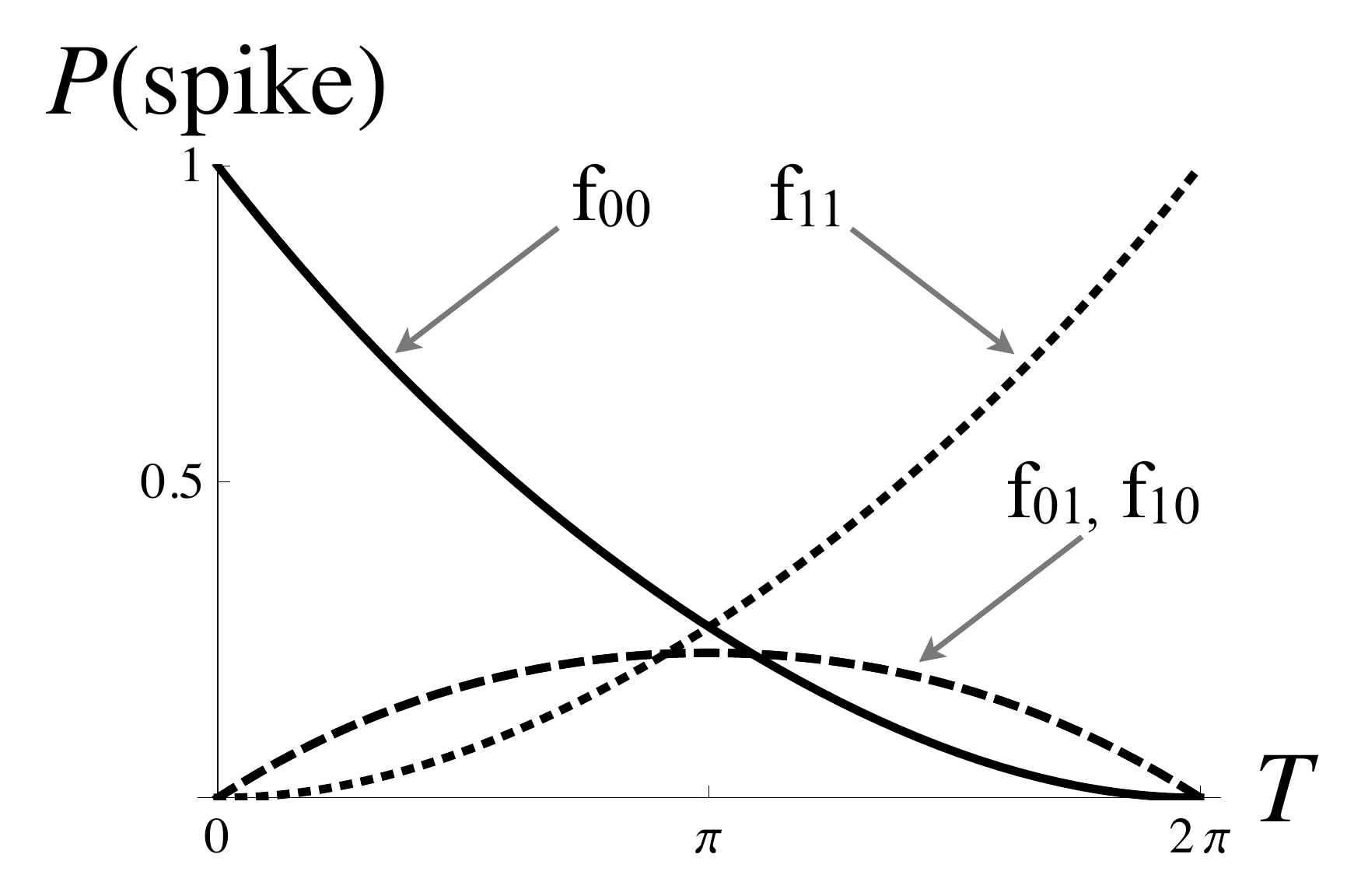}
\caption{Joint spiking probability for two oscillators receiving partially correlated noise is shown for observations windows $T \leq 1$, where $1$ is the natural frequency of the oscillation. The subscripts $ij$ indicate the probability that oscillator $i$ or $j$ does (1) or does not (0) spike.}
\label{fig:probabilitySpike}
\end{center}
\end{figure}
%------------------------------------

Now we will calculate the spike count correlation directly for observation windows $T$ that are shorter than or equal to the natural period, which we will assume to be $2\pi$. First let us consider the probability that a spike occurs in $[0,T]$. We say that oscillator $i$ spikes when its phase $\theta_i$ reaches $2\pi$, which is to say $\theta_i(T) \geq 2\pi$. Assuming as usual that the noise amplitude $\sigma$ is small, we expand the phase to lowest order as in Eq.(\ref{eq:lowestOrder}), that is $\theta_i(T) = \theta_i(0) + T + \mathcal{O}(\sigma)$. Therefore the probability that oscillator $i$ spikes is simply
\begin{align*}
\text{P}[\theta_i \text{ spikes}] &= \text{P}[\theta_i + T \geq 2\pi]\\
\text{P}[\theta_i \text{ does not spike}] &= \text{P}[\theta_i + T < 2\pi].
\end{align*} 
For two oscillators, there are four possibilities for the joint spike count:
%\begin{widetext}
\begin{align*}
\text{P}[\theta_1& \text{ does not spike},  \theta_2 \text{ does not spike}] \\
&\indent = \text{P}[\theta_1 + T < 2\pi,  \theta_2 + T < 2\pi]\\
 \text{P}[\theta_1& \text{ spikes},  \theta_2 \text{ does not spike}] \\
&\indent = \text{P}[\theta_1 + T \geq 2\pi,  \theta_2 + T < 2\pi]\\
\text{P}[\theta_1& \text{ does not spike},  \theta_2 \text{ spikes}] \\
&\indent = \text{P}[\theta_1 + T < 2\pi,  \theta_2 + T \geq 2\pi]\\
\text{P}[\theta_1& \text{ spikes},  \theta_2 \text{ spikes}] \\
&\indent = \text{P}[\theta_1 + T \geq 2\pi,  \theta_2 + T \geq 2\pi].
\end{align*} 
%\end{widetext}
%------------------------------------
\begin{figure*}
\begin{center}
\includegraphics[width=6.5in]{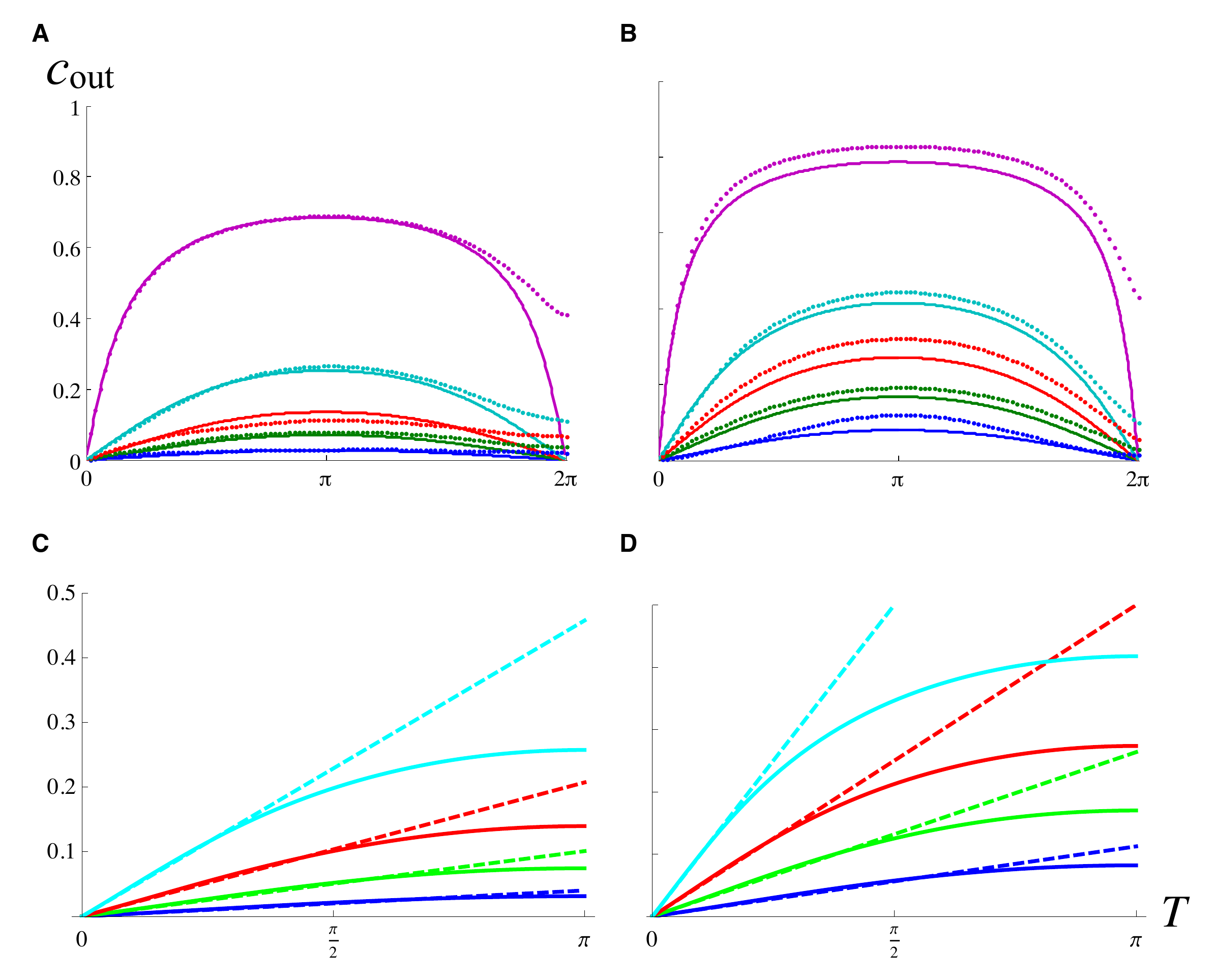}
\caption{(A,B) Theoretical (solid) and simulated (dotted) output correlation curves are shown as a function of the observation window $T\leq 2\pi$. (A) Type I oscillators. (B) Type II oscillators. (C,D) The initial slope (dashed) of the spike count correlation (solid) is the linear approximation of Eq.(\ref{eq:smallTcorr}) at $T=0$, which is given in Eq.(\ref{eq:susceptibility}). (C) Type I oscillators. (D) Type II oscillators. For all plots, noise amplitude $\sigma=0.04$, and colors indicate the level of input correlation: 0.2 (blue), 0.4 (green), 0.6 (red), 0.8 (cyan), 0.99 (purple). In all cases, noise amplitude $\sigma = 0.05$.}
\label{fig:small_T}
\end{center}
\end{figure*}
%------------------------------------

These probabilities can be obtained directly by integrating the density of the phase difference, Eq.(\ref{eq:ssPhaseDistribution}), over the appropriate domain. Note that this gives four discrete joint probabilities for each observation window $T\in [0,2\pi]$. For convenience, let us define the following functions of $T$:
\begin{align*}
f_{00}(T) &:= \text{P}[\theta_1 \leq 2\pi - T,  \theta_2 \leq 2\pi  - T]\\
 &= \frac{1}{2\pi}\int_0^{2\pi  - T} \int_0^{2\pi - T} P(y-x) dx dy \\
f_{01}(T) &:= \text{P}[\theta_1 > 2\pi - T,  \theta_2 \leq 2\pi - T] \\
&= \frac{1}{2\pi}\int_{2\pi  - T}^{2\pi} \int_0^{2\pi - T} P(y-x) dx dy \\
f_{10}(T) &:= \text{P}[\theta_1 \leq 2\pi - T,  \theta_2 > 2\pi - T] \\
&= \frac{1}{2\pi}\int_0^{2\pi  - T} \int_{2\pi - T}^{2\pi} P(y-x) dx dy \\
f_{11}(T) &:= \text{P}[\theta_1 > 2\pi - T,  \theta_2 > 2\pi - T] \\
&= \frac{1}{2\pi}\int_{2\pi  - T}^{2\pi} \int_{2\pi  - T}^{2\pi} P(y-x) dx dy.
\end{align*} 

Let $X$ be the random variable such that $X=1$ if $\theta_1$ spikes during the observation period $T$, and $X=0$ if $\theta_1$ does not spike. Similarly, let $Y$ represent the presence or absence of a spike in oscillator~$\theta_2$. Then the covariance is given by $\text{Cov}[X,Y] = \text{E}[XY] -\text{E}[X]\text{E}[Y]$. In terms of the functions defined above we have
\begin{align*}
\text{E}[X] &= 0 \cdot (f_{00}+ f_{01}) + 1 \cdot (f_{10}+ f_{11}) \\
&= (f_{10}+ f_{11}) = \text{E}[X^2] \\
\text{E}[Y] &= 0 \cdot (f_{00} + f_{10}) + 1 \cdot (f_{01} + f_{11}) \\
&= (f_{01} + f_{11}) = \text{E}[Y^2] \\
\text{E}[XY] &= 0 \cdot 0 \cdot f_{00} + 1 \cdot 0 \cdot f_{10} + 0 \cdot 1 \cdot f_{01} + 1 \cdot 1 \cdot f_{11} \\
&= f_{11}.
\end{align*} 

A few simplifications are possible. In particular, the sum $f_{10}(T)+ f_{11}(T)$ is just the marginal probability that $\theta_1$ spikes within time $T$. Since $\theta_1$ is uniformly distributed, this probability is simply $\frac{T}{2\pi}$. Furthermore, we also have $f_{10}=f_{01}$ by the symmetry of the density $P$, and hence $\sqrt{\text{Var} [X] \text{Var} [Y]} = \text{Var} [X]$. Therefore  the spike count correlation over short time windows is 

\begin{align}
& \text{Cor}[X,Y](T;c)  \\
&= \frac{\text{E}[XY] -\text{E}[X]\text{E}[Y]}{\text{Var} [X]} \nonumber\\
&= \frac{f_{11} - (f_{10}+ f_{11})^2}{(f_{10}+ f_{11})(1-(f_{10}+ f_{11}))} \nonumber\\
&= \frac{f_{11} - \left(\frac{T}{2\pi}\right)^2}{\frac{T}{2\pi} \left(1-\frac{T}{2\pi}\right)} \nonumber\\
&= \frac{1}{2\pi T-T^2} \left[ 2\pi \int_{2\pi-T}^{2\pi}\int_{2\pi-T}^{2\pi} P(y-x)dx dy - T^2\right].
\label{eq:smallTcorr}
\end{align}

Fig.(\ref{fig:small_T}A,B) shows how this analytically derived output correlation compares with numerical simulations for type I and type II oscillators, respectively. 

We can make a further simplification by considering the linear part of Eq.(\ref{eq:smallTcorr}) for $T$ close to zero:

\begin{equation*}
c_{out} = T \left(P(0) - \frac{1}{2\pi}\right) + \mathcal{O}(T^2)
\end{equation*}
Thus, the initial slope of the output correlation is proportional to the peak of the stationary phase difference distribution, $P(\phi)|_{\phi=0}$. Substituting $P_I(0)$ and $P_{II}(0)$ from Eq.(\ref{eq:PI}) and Eq.(\ref{eq:PII}), we obtain:

\begin{align}
c_{out,I}  &\approx \frac{T}{\pi} \left(  \frac{c}{3(1-c)+\sqrt{3(c-1)(c-3)}} \right) \nonumber\\
 &= T\frac{c}{6\pi} + \mathcal{O}(c^2) \nonumber\\
c_{out,II}  &\approx  \frac{T}{2\pi} \left( \frac{1+c}{\sqrt{1-c^2}} -1 \right) \\
 &= T\frac{c}{2\pi} + \mathcal{O}(c^2).
\label{eq:susceptibility}
\end{align}
From here, it is clear that the initial slope of $c_{out}$ is greater for type II than for type I oscillators; in fact the type II output correlation rises three times faster than the type I, to lowest order in $c$. See Fig.(\ref{fig:small_T}C,D).

%------------------------------------
\begin{figure}
\begin{center}
\includegraphics[width=3.25in]{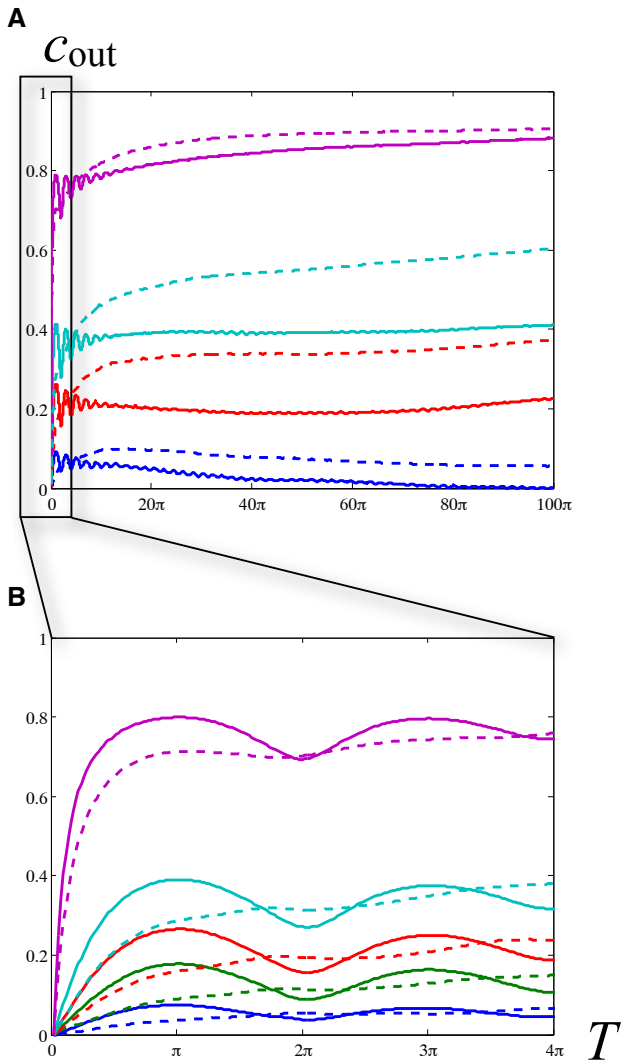}
\caption{Output correlation is shown as a function of intermediate-length observation windows $T$. Colors indicate the level of input correlation: 0.2 (blue), 0.4 (green), 0.6 (red), 0.8 (cyan), 0.99 (purple). (A) Type II oscillators (solid) exhibit higher output correlations over short time scales than do type I (dashed). (B) This result reverses over long time scales. In all cases, noise amplitude $\sigma = 0.2$.}
\label{fig:sim_type1&2}
\end{center}
\end{figure}
%------------------------------------

%%%%%%%%%%%%%%%%%%%%%%%%%%%%%%%%%%%%%%%%%%%%%%%%%%%%
\section*{Discussion}
%%%%%%%%%%%%%%%%%%%%%%%%%%%%%%%%%%%%%%%%%%%%%%%%%%%%

We have demonstrated a novel approach to approximating the spike count correlation of noisy neural oscillators over both long and short time scales. In the case of long windows of observation $T$ much greater than the natural period of oscillation, we used the total elapsed phase, modulo the period, as a proxy for the spike count. The difference between these quantities is at most one and hence is negligible for when many spikes are observed over large time windows $T$. In our perturbation expansion to lowest order in the noise amplitude, $\sigma$, the correlation between oscillators depends only on the PRC and the stationary distribution of the phase difference. A further approximation assuming small input correlation $c$ reveals that output correlation scales with the autocorrelation of the PRC, which is a nonnegative quantity that equals zero precisely when the PRC is a pure sinusoid, i.e., when the oscillator displays type II dynamics. This observation sheds some light on the counterintuitive finding, first reported by Barreiro,  et al. \cite{barreiro:2010}, whereby type I oscillators transfer correlations more faithfully than do type II over long time scales, although the reverse holds true for the better understood case of short time scales.

Using straightforward probabilistic reasoning, we computed the spike count correlation directly for short time scales. In the limit of small $T$ and small $c$, we obtain an expression for the initial slope of the output correlation, also known as the correlation susceptibility \cite{brent:nature}. In \cite{brent:nature}, de la Rocha, et al. use a phenomenological model to explore the complex relationship between susceptibility, firing rate and threshold nonlinearities. The present analysis illustrates the contribution of bifurcation structure via phase resetting dynamics. In particular, the susceptibility is proportional to the peak of the stationary phase difference distribution, $P(\phi)|_{\phi=0}$, which in turn depends on the shape of the PRC.

Our analytic expressions in the limit of small noise agree well with spike count correlations computed from simulated oscillators. However, for tractability we included only terms of order one in the perturbation expansion of the phase given in Eq.(\ref{eq:lowestOrder}). As a result, the present analysis cannot account for the slow drift of the correlation due to noise, which is visible for values of $T$ near $2\pi$ in Fig.(\ref{fig:small_T}), and is even more apparent for the intermediate values of $T$ shown in Fig.(\ref{fig:sim_type1&2}). New analytic methods capable of addressing non-extremal cases would shed light on this and many other questions in mathematical biology.
\vskip 0.01mm
\bibliography{spk_cnt_paper}
\end{document}